\documentclass[10pt,a4paper]{amsart}

\usepackage{graphicx,amssymb,color,amscd}
\usepackage[all]{xy}

\theoremstyle{plain}
\newtheorem{theorem}{Theorem}[section]
\newtheorem{lemma}[theorem]{Lemma}
\newtheorem{proposition}[theorem]{Proposition}
\newtheorem{corollary}[theorem]{Corollary}
\theoremstyle{definition}
\newtheorem{definition}[theorem]{Definition}
\theoremstyle{remark}
\newtheorem{remark}[theorem]{Remark}
\newtheorem*{acknowledgments}{Acknowledgments}
\numberwithin{equation}{section}
\numberwithin{figure}{section}

\newcommand\nc{\newcommand}
\nc\rnc{\renewcommand}

\nc\Br{\operatorname{Br}}
\nc\BrJ{\Br(\bJ_{2n}(n+1))}

\nc\modZ {\mathbb{Z}}
\nc\modN {\mathbb{N}}
\nc\modQ {\mathbb{Q}}
\nc\modL {\mathcal{L}}
\nc\zzzcolon {\colon}
\nc\zzzvert {\;|\;}
\nc\bB{\bar{B}}
\nc\bJ{\bar{J}}
\nc\bJB{\bJ^B}
\nc\Span{\operatorname{Span}}
\nc\modS {{\mathcal S}}
\nc\modF {{\mathcal F}}
\nc\modT {{\mathcal T}}
\nc\modG {{\mathcal G}}
\nc\modB {{\mathcal B}}
\nc\modP {{\mathcal P}}
\nc\modM {{\mathcal M}}
\nc\modone {{\mathbf 1}}
\nc\FI[2]{\begin{figure}
    \begin{center}\input{#1.pstex_t}\end{center}
    \caption{#2}
    \label{#1}
  \end{figure}}
\nc\xto[1]{\overset{#1}{\longrightarrow}}
\nc\xyto[2]{\overset{#1}{\underset{#2}{\longrightarrow}}}
\nc\simeqto{\xto{\simeq}}
\nc\half{\frac12}
\nc\Sym{\operatorname{Sym}}
\nc\SymT{\tilde\Sym^2\modT _{n+1}}
\nc\sums{\sum_{\sigma \in S_{n-1}}}
\nc\sumss{\sum_{\sigma ,\sigma '\in S_{n-1}}}
\nc\SLnn{SL_{n+1}}
\nc\BSLnn{BSL_{n+1}}
\nc\bBSLnn{\overline{BSL}_{n+1}}
\nc\Bnn{B_{n+1}}
\nc\bBnn{\overline{B}_{n+1}}

\begin{document}

\title[Finite type invariants and Milnor invariants for Brunnian
  links]{Finite type invariants and Milnor invariants\\for Brunnian
  links}

\author[K. Habiro]{Kazuo Habiro}

\address{Research Institute for Mathematical Sciences\\Kyoto
  University\\Kyoto 606-8502, Japan}

\email{habiro@kurims.kyoto-u.ac.jp}

\author[J.-B. Meilhan]{Jean-Baptiste Meilhan}

\email{meilhan@kurims.kyoto-u.ac.jp}

\thanks{The first author is partially supported by the Japan Society
  for the Promotion of Science, Grant-in-Aid for Young Scientists (B),
  16740033.  The second author is supported by a Postdoctoral
  Fellowship and a Grant-in-Aid for Scientific Research of the Japan
  Society for the Promotion of Science.}

\date{April 3, 2007 (First version: October 25, 2005)}
%
%

\begin{abstract}
  A link $L$ in the $3$-sphere is called {\em Brunnian} if every
  proper sublink of $L$ is trivial.  In a previous paper, the first
  author proved that the restriction to Brunnian links of any
  Goussarov-Vassiliev finite type invariant of $(n+1)$-component
  links of degree $<2n$ is trivial.  The purpose of this paper is to
  study the first nontrivial case.  We show that the restriction
  of an invariant of degree $2n$ to $(n+1)$-component Brunnian links
  can be expressed as a quadratic form on the Milnor link-homotopy
  invariants of length $n+1$.
\end{abstract}

\maketitle

\section{Introduction}
\label{sec:introduction}

The notion of {\em Goussarov-Vassiliev finite type link invariants}
\cite{Gusarov:91,Gusarov:94,Vassiliev} enables us to understand
the various quantum invariants from a unifying viewpoint, see
e.g. \cite{BN,O2}.  The theory involves a descending filtration
\begin{equation*}
  \modZ \modL (m)=J_0(m)\supset J_1(m)\supset \dots
\end{equation*}
of the free abelian group $\modZ \modL (m)$ generated by the set $\modL (m)$ of
the ambient isotopy classes of $m$-component, oriented, ordered links
in $S^3$.  Here each $J_n(m)$ is generated by alternating sums of
links over $n$ independent crossing changes.  A homomorphism from
$\modZ \modL (m)$ to an abelian group $A$ is said to be a Goussarov-Vassiliev
invariant of degree $n$ if it vanishes on $J_{n+1}(m)$.  Thus, for
$L,L'\in \modL (m)$, we have $L-L'\in J_{n+1}(m)$ if and only if $L$ and $L'$
have the same values of Goussarov-Vassiliev invariants of degree $\le n$
with values in any abelian group.

It is natural to ask what kind of informations a Goussarov-Vassiliev
link invariants can contain and what is the topological meaning of the
unitrivalent diagrams.  {\em Calculus of claspers}, introduced by
Goussarov and the first author \cite{Goussarov99,Goussarov01,H},
answers these questions.  (We will recall the definition of claspers
in Section \ref{sec:claspers}.)  A special type of claspers, called
{\em graph claspers}, can be regarded as topological realizations of
unitrivalent diagrams.  For {\em knots}, claspers enables us to give a
complete topological characterization of the informations that can be
contained by Goussarov-Vassiliev invariants of degree $<n$
\cite{Goussarov01,H}: The difference of two knots is in $J_n$ if and
only if these two knots are {\em $C_n$-equivalent}.  Here
$C_n$-equivalence is generated by a certain type of local moves,
called {\em $C_n$-moves} (called {\em $(n-1)$-variations} by
Goussarov), which is defined as surgeries along certain tree claspers.

For links with more than $1$ components, the above-mentioned
properties of Goussarov-Vassiliev invariants does not hold.  It is
true that if $L,L'\in \modL (m)$ are $C_n$-equivalent, then we have
$L-L'\in J_n(m)$, but the converse does not hold in general.  A
counterexample is Milnor's link $L_{n+1}$ of $n+1$ components depicted
in Figure \ref{F34}: If $n\ge 2$, $L_n$ is ($C_n$-equivalent but) not
$C_{n+1}$-equivalent to the $(n+1)$-component unlink $U$, while we
have $L_{n+1}-U\in J_{2n}(n+1)$ (but $L_{n+1}-U\not\in J_{2n+1}(n+1)$),
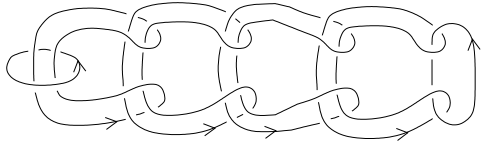
see \cite[Proposition 7.4]{H}.  \FI{F34}{Milnor's link $L_6$ of $6$
components} (This fact is contrasting to the case of {\em string
links}: Conjecturally \cite[Conjecture 6.13]{H}, two string links
$L,L'$ of the same number of components are $C_n$-equivalent if and
only if $L-L'\in J_n$.)

Milnor's links are typical examples of {\em Brunnian links}.  Recall
that a link in an oriented, connected $3$-manifold is said to be
Brunnian if every proper sublink of it is an unlink.  In some sense,
an $n$-component Brunnian link is a `pure $n$-component linking'.
Thus studying the behavior of Goussarov-Vassiliev invariants on
Brunnian links would be a first step in understanding the
Goussarov-Vassiliev invariants for links.

The first author generalized a part of the above-mentioned properties
of Milnor's links to Brunnian links:

\begin{theorem}[\cite{Hb}]
  \label{r1}
  Let $L$ be an $(n+1)$-component Brunnian link in a connected,
  oriented $3$-manifold $M$ ($n\ge 1$), and let $U$ be an
  $(n+1)$-component unlink in $M$.  Then we have the following.
  \begin{enumerate}
  \item $L$ and $U$ are $C_n$-equivalent.
  \item If $n\ge 2$, then we have $L-U\in J_{2n}(n+1)$.  Hence $L$ and $U$
  are not distinguished by any Goussarov-Vassiliev invariants of
  degree $<2n$.
  \end{enumerate}
\end{theorem}

The case $M=S^3$ of Theorem \ref{r1} was announced in \cite{H}, and
was later proved also by Miyazawa and Yasuhara \cite{MY},
independently to \cite{Hb}.

The purpose of the present paper is to study the restrictions of
Goussarov-Vassiliev invariants of degree $2n$ to $(n+1)$-component
Brunnian links in $S^3$, which is the first nontrivial case according
to Theorem \ref{r1}.  The main result in the present paper expresses
any such restriction as a {\em quadratic} form of Milnor $\bar\mu$ link-homotopy
invariants of length $n+1$:

\begin{theorem}
  \label{r6}
  Let $f$ be any $\modZ $-valued Goussarov-Vassiliev link invariant of
  degree $2n$.  Then there are (non-unique) integers $f_{\sigma ,\sigma '}$ for
  elements $\sigma ,\sigma '$ of the symmetric group $S_{n-1}$ on the set
  $\{1,\ldots ,n-1\}$ such that, for any $(n+1)$-component Brunnian link
  $L$, we have
  \begin{equation}
    \label{e8}
    f(L)-f(U)= \sumss f_{\sigma ,\sigma '}\bar\mu _\sigma (L)\bar\mu _{\sigma '}(L).
  \end{equation}
  Here, $U$ is an $(n+1)$-component unlink, and we
  set, for $\sigma \in S_{n-1}$,
  \begin{equation*}
    \bar\mu _\sigma (L)= \bar\mu _{\sigma (1),\sigma (2),\ldots ,\sigma (n-1),n,n+1}(L)\in \modZ.
  \end{equation*}
\end{theorem}

A possible choice for the integers $f_{\sigma ,\sigma '}$ is given in terms of tree claspers
in Section \ref{sec:proof-theorem-refr6}.

\begin{remark}
The proof of Theorem \ref{r6} involves calculus of claspers.
The first preprint version of the present paper
(arxiv:math.GT/0510534v1) contained a one-page sketch of an
alternative proof of Theorem \ref{r6} using the Kontsevich integral.  This
alternative proof has been separated from the present paper, and has
been published in \cite{HM2}.  
Though shorter than the clasper-based proof below, the proof in \cite{HM2} relies heavily on the properties of
unitrivalent diagrams, and the topological meaning of the steps in the
proof are therefore not always very clear.  
The present proof gives a better understanding of Theorem \ref{r6} from that point of view.  
\end{remark}

Recall that Milnor invariants of length $n+1$ for string links are
Goussarov-Vassiliev invariants of degree $\le n$ \cite{BN2,Lin} (see
also \cite{HM}).  As is well-known, Milnor's invariants is not
well-defined for all links, and hence it does not make sense to ask
whether Milnor invariants of length $n+1$ is of degree $\le n$ or not.
However, as Theorem \ref{r6} indicates, a quadratic expression in such
Milnor invariants, which is well-defined at least for
$(n+1)$-component Brunnian links, may extend to a link invariant of
degree $\le 2n$.

In the study of Milnor's invariants, tree claspers seem at least as
useful as Cochran's construction \cite{Cochran}.  For the use of
claspers in the study of the Milnor invariants, see also
\cite{Gar1,Habe,Meilhan}.  For other relationships between finite type
invariants and the Milnor invariants, see
\cite{BN2,Lin,HM,HO2,Masbaum-Vaintrob}.

\vskip 5pt

We organize the rest of the paper as follows.

In Section \ref{sec:claspers}, we recall some definitions from clasper
calculus.

In Section \ref{sec:cak-equivalence}, we recall the notion of
$C^a_k$-equivalence for links, studied in \cite{Hb}.  If a link $L$ is
$C^a_k$-equivalent (for any $k$) to a Brunnian link, then $L$ also is
a Brunnian link.

In Section \ref{sec:n+1-comp-brunn}, we study the group $\bBSLnn$ of
$C^a_{n+1}$-equivalence classes of $(n+1)$-component string links.  We
establish an isomorphism
\begin{equation*}
  \theta _n\zzzcolon \modT _{n+1}\simeqto\bBSLnn
\end{equation*}
from an abelian group $\modT _{n+1}$ of certain tree diagrams.  This map
is essentially the inverse to the Milnor link-homotopy invariants of
length $n+1$.

In Section \ref{sec:group-bbnn}, we apply the results in Section
\ref{sec:n+1-comp-brunn} to Brunnian links.  The operation of closing
string links induces a bijection
\begin{equation*}
  \bar{c}_{n+1}\zzzcolon \bBSLnn\simeqto\bBnn,
\end{equation*}
where $\bBnn$ is the set of $C^a_{n+1}$-equivalence classes of
$(n+1)$-component Brunnian links.  As a byproduct, we obtain another
proof of a result of Miyazawa and Yasuhara \cite{MY}.

In Section \ref{sec:gouss-vass-filtr}, we recall the definition of the
Goussarov-Vassiliev filtration for links using claspers.

In Section \ref{sec:gouss-vass-invar}, we study the behavior of
Goussarov-Vassiliev invariants of degree $2n$ for $(n+1)$-component
Brunnian links.  We first show that two $C^a_{n+1}$-equivalent,
$(n+1)$-component Brunnian links cannot be distinguished by
Goussarov-Vassiliev invariants of degree $2n$.  We have a quadratic map
\begin{equation*}
  \kappa _{n+1}: \bBnn \longrightarrow \bJ_{2n}(n+1)
\end{equation*}
defined by $\kappa _{n+1}([L]_{C^a_{n+1}})=[L-U]_{J_{2n+1}}$.  We prove
Theorem \ref{r6}, using $\kappa _{n+1}$.

\begin{acknowledgments}
  The authors wish to thank Akira Yasuhara
  for helpful conversations.
\end{acknowledgments}

\section{Claspers}
\label{sec:claspers}

In this section, we recall some definitions from calculus of claspers.
For the details, we refer the reader to \cite{H}.

A {\em clasper} in an oriented $3$-manifold $M$ is a compact, possibly
unorientable, embedded surface $G$ in $\operatorname{int} M$ equipped with a
decomposition into connected subsurfaces called {\em leaves}, {\em
disk-leaves}, {\em nodes}, {\em boxes}, and {\em edges}.  Two distinct
non-edge subsurfaces are disjoint.  Edges are disjoint bands which
connect two subsurfaces of the other types.  A connected component of
the intersection of one edge $E$ and another subsurface $F$ (of
different type), which is an arc in $\partial E\cap \partial F$, is called an {\em
attaching region} of $F$.
\begin{itemize}
\item A {\em leaf} is an annulus with one attaching region.
\item A {\em disk-leaf} is a disk with one attaching region.
\item A {\em node} is a disk with three attaching regions.
  (Usually, a node is incident to three edges, but it is allowed
  that the two ends of one edge are attached to a node.)
\item A {\em box} is a disk with three attaching regions.  (The same
  remark as that for node applies here, too.)  Moreover, one attaching
  region is distinguished with the other two.  (This distinction is
  done by drawing a box as a rectangle, see \cite{H}.)
\end{itemize}
A clasper $G$ {\em for} a link $L$ in $M$ is a clasper in $M$ such
that the intersection $G\cap L$ consists of finitely many transverse
double points and is contained in the interior of the union of
disk-leaves.

We often use the drawing convention for claspers as described in
\cite{H}.

{\em Surgery along a clasper} $G$ is defined to be surgery along the
{\em associated framed link $L_G$} to $G$.  Here $L_G$ is obtained
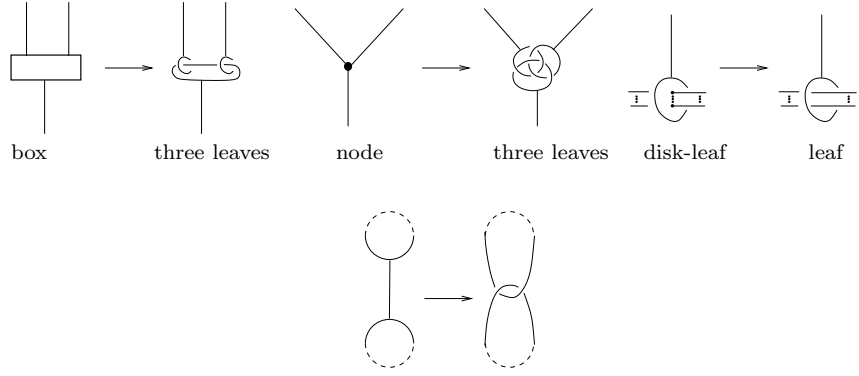
from $G$ by the rules described in Figure \ref{F27new}.  \FI{F27new}{How to
obtain the associated framed link $L_G$ from $G$.  First one replaces
boxes, nodes and disk-leaves by leaves.  Then replace each
`I-shaped' clasper by a $2$-component framed link as depicted.}


A {\em tree clasper} is a connected clasper $T$ without boxes, such
that the union of edges and nodes of $T$ is simply connected.  A tree
clasper $T$ is called {\em strict} if each component of $T$ has no
leaves and at least one disk-leaf.  Surgery along a strict tree
clasper $T$ is {\em tame} in the sense of \cite[Section 2.3]{H}, i.e.,
the result of surgery along $T$ preserves the $3$-manifold and the
surgery may be regarded as a move on a link.

A tree clasper $T$ for a link $L$ is {\em simple} (with respect to
$L$) if each disk-leaf of $T$ has exactly one intersection point with
$L$.

The {\em degree} of a strict tree clasper $G$ is defined to be the
number of nodes of $T$ plus $1$.  For $n\ge 1$, a {\em $C_n$-tree} is
a strict tree clasper of degree $n$.  A {\em (simple) $C_n$-move} is a
local move on links defined as surgery along a (simple) $C_n$-tree.
For example, a simple $C_1$-move is a crossing change, and a simple
$C_2$-move is a delta move \cite{Matveev,Murakami-Nakanishi}.  The
{\em $C_n$-equivalence} is the equivalence relation on links generated
by $C_n$-moves.  This equivalence relation is also generated by simple
$C_n$-moves.  The $C_n$-equivalence becomes finer as $n$ increases.

\section{$C^a_k$-equivalence}
\label{sec:cak-equivalence}

We recall from \cite{Hb} the definition of the $C^a_k$-equivalence.

\begin{definition}
  \label{r9}
  Let $L$ be an $m$-component link in a $3$-manifold $M$.  For $k\ge m-1$,
  a {\em $C^a_k$-tree} for $L$ in $M$ is a $C_k$-tree $T$ for $L$ in
  $M$, such that
  \begin{enumerate}
  \item for each disk-leaf $A$ of $T$, all the strands intersecting
    $A$ are contained in one component of $L$, and
  \item each component of $L$ intersects at least one disk-leaf of
  $T$, i.e., $T$ intersects all the components of $L$.
  \end{enumerate}
  Note that the condition (1) is vacuous if $T$ is simple.

  A {\em $C^a_k$-move} on a link is surgery along a $C^a_k$-tree.  The
  {\em $C^a_k$-equivalence} is the equivalence relation on links
  generated by $C^a_k$-moves.  A {\em $C^a_k$-forest} is a clasper
  consisting only of $C^a_k$-trees.
\end{definition}

Clearly, the above notions are defined also for tangles, particularly
for string links.

What makes the notion of $C^a_k$-equivalence useful in the study of
Brunnian links is the fact that a link which is $C^a_k$-equivalent
(for any $k$) to a Brunnian link is again a Brunnian link
(\cite[Proposition 5]{Hb}).

Note that the $C^a_k$-equivalence is generated by {\em simple}
$C^a_k$-moves, i.e., surgeries along simple $C^a_k$-trees \cite{Hb}.
In the following, we use technical lemmas from \cite{Hb}.

\begin{lemma}[{\cite[Lemma 7]{Hb}, $C^a$-version of \cite[Theorem
	3.17]{H}}]
  \label{r18}
  For two tangles $\beta $ and $\beta '$ in a $3$-manifold $M$, and an integer
  $k\ge 1$, the following conditions are equivalent.
  \begin{enumerate}
  \item $\beta $ and $\beta '$ are $C^a_k$-equivalent.
  \item There is a simple $C^a_k$-forest $F$ for $\beta $ in $M$ such
    that $\beta _F\cong \beta '$.
  \end{enumerate}
\end{lemma}

\begin{lemma}[{\cite[Lemma 8]{Hb}, $C^a$-version of \cite[Proposition
    4.5]{H}}]
  \label{r29}
  Let $\beta $ be a tangle in a $3$-manifold $M$, and let $\beta _0$ be a
  component of $\beta $.  Let $T_1$ and $T_2$ be $C_k$-trees for a tangle
  $\beta $ in $M$, differing from each other by a crossing change of an
  edge with the component $\beta _0$.  Suppose that $T_1$ and $T_2$ are
  $C^a_k$-trees for either $\beta $ or $\beta \setminus \beta _0$.  Then $\beta _{T_1}$ and
  $\beta _{T_2}$ are related by one $C^a_{k+1}$-move.
\end{lemma}

\section{The group $\bBSLnn$}
\label{sec:n+1-comp-brunn}

\subsection{The monoids $\BSLnn$ and $\bBSLnn$}
\label{sec:groups-bbn+1}

Let us recall the definition of string links.  (For the details, see
e.g. \cite{Habegger-Lin,H}).  Let $x_1,\ldots ,x_{n+1}\in \operatorname{int} D^2$
be distinct points.  An $(n+1)$-component {\em string link}
$\beta =\beta _1\cup \dots \cup \beta _{n+1}$ is a tangle in the cylinder $D^2\times [0,1]$,
consisting of arc components $\beta _1,\ldots ,\beta _{n+1}$ such that
$\partial \beta _i=\{x_i\}\times \{0,1\}$ for each $i$.  Let $\SLnn$ denote the set of
$(n+1)$-component string links up to ambient isotopy fixing endpoints.
There is a natural, well-known monoid structure for $\SLnn$ with
multiplication given by `stacking' of string links.
The identity string link is denoted by $\modone =\modone _{n+1}$.

Let $\BSLnn$ denote the submonoid of $\SLnn$ consisting of {\em
  Brunnian string links}.  Here a string link $\beta $ is said to be
Brunnian if every proper subtangle of $\beta $ is the identity string
link.

We have the following characterization of Brunnian string links.

\begin{theorem}[{\cite[Theorem 9]{Hb}, \cite[Proposition 4.1]{MY}}]
  \label{r12}
  An $(n+1)$-component link (resp. string link) is Brunnian if and
  only if it is $C^a_n$-trivial, i.e., it is $C^a_n$-equivalent to
  the unlink (resp. the identity string link).
\end{theorem}

Set
\begin{equation*}
  \bBSLnn=\BSLnn/(\text{$C^a_{n+1}$-equivalence}).
\end{equation*}
By Theorem \ref{r12}, $\bBSLnn$ can be regarded as the monoid of
$C^a_{n+1}$-equivalence classes of $C^a_n$-trivial, $(n+1)$-component
string links (in $D^2\times [0,1]$).

In the rest of this section, we will describe the structure of
$\bBSLnn$.

\subsection{The group $\bBSLnn$ and the surgery map $\theta _n\zzzcolon \modT _{n+1}\rightarrow \bBSLnn$}
\label{sec:groups-bbslnns}

\begin{proposition}
  \label{r11}
  $\bBSLnn$ is a finitely generated abelian group.
\end{proposition}

\begin{proof}
  The assertion is obtained by adapting the proof of \cite[Lemma 5.5,
  Corollary 5.6]{H} into the $C^a$ setting.
\end{proof}

Let $n\ge 1$.  By a {\em (labeled) unitrivalent tree} of degree $n$ we
mean a vertex-oriented, unitrivalent graph $t$ such that the $n+1$
univalent vertices of $t$ are labeled by distinct elements from
$\{1,2,\ldots ,n+1\}$.  In figures, the counterclockwise
vertex-orientation is assumed at each vertex.

Let $\modT _{n+1}$ denote the free abelian group generated by
unitrivalent trees of degree $n$, modulo the well-known IHX and AS
relations.

For a unitrivalent tree $t$, let $T_t$ denote a $C^a_n$-tree for $\modone $
such that the tree shape and the labeling of $T_t$ is induced by those
of $t$, and such that after choosing an orientation of $T_t$, for each
$i=1,\ldots ,n+1$, the sign of the intersection of the $i$th string of $\modone $
and the disk-leaf of $T_t$ corresponding to the univalent vertex of
$t$ colored $i$ is positive.  See for example Figure \ref{F35}.
\FI{F35}{}

\begin{proposition}
  \label{r25}
  There is a unique isomorphism
  \begin{equation*}
    \theta _{n+1}\zzzcolon \modT _{n+1}\simeqto\bBSLnn.
  \end{equation*}
  such that $\theta _{n+1}(t)=[\modone _{T_t}]_{C^a_{n+1}}$ for each unitrivalent
  tree $t$, where $T_t$ is as above.
\end{proposition}

\begin{proof}
  Let $\modT '_{n+1}$ be the free abelian group generated by unitrivalent
  trees of degree $n$, modulo the AS relations.  By adapting the proof
  of \cite[Theorem 4.7]{H} into the $C^a$ setting, we see that there is
  a unique surjective homomorphism
  \begin{equation*}
    \theta '_{n+1}\zzzcolon \modT '_{n+1}\rightarrow \bBSLnn.
  \end{equation*}

  To see that $\theta '_{n+1}$ factors through the projection
  $\modT'_{n+1}\rightarrow \modT _n$, it suffices to see that the IHX
  relation is valid in $\bBSLnn$, i.e., $t_I-t_H+t_X\in \modT '_{n+1}$
  is mapped to $0$, where $t_I,t_H,t_X$ locally differs as in the
  definition of the IHX relation.  This can be checked by adapting the
  IHX relation for tree claspers (see e.g. \cite{Goussarov01,CT,CST})
  into the $C^a$ setting.

  Let
  \begin{equation*}
    \theta _{n+1}\zzzcolon \modT _{n+1}\rightarrow \bBSLnn
  \end{equation*}
  be the surjective homomorphism induced by $\theta '_{n+1}$.  As in the
  statement of Theorem \ref{r6}, for $\sigma \in S_{n-1}$ and $L\in B(n+1)$, we
  set
  \begin{equation*}
    \mu _\sigma (T) = \mu _{\sigma (1),\sigma (2),\ldots ,\sigma (n-1),n,n+1}(T),
  \end{equation*}
  where $\mu _{\sigma (1),\sigma (2),\ldots ,\sigma (n-1),n,n+1}(T)\in \modZ $ is the Milnor
  string link invariant of $T$.  Let $t_\sigma $ denote the unitrivalent
  tree as depicted in Figure \ref{F32}.
  {\def\te#1{$\sigma (#1)$}\FI{F32}{}}  The $t_\sigma $ for $\sigma \in S_{n-1}$ form a
  basis of $\modT _{n+1}$.
  Define a homomorphism
  \begin{equation*}
    \mu _{n+1}\zzzcolon \bBSLnn\longrightarrow\modT _{n+1}
  \end{equation*}
  by
  \begin{equation*}
    \mu _{n+1}(L) = \sum_{\sigma \in S_{n-1}} \mu _\sigma (L) t_\sigma ,
  \end{equation*}
  By \cite[Theorem 7.2]{H}, $\mu _{n+1}$ is well defined.

  To show that $\mu _{n+1}$ is left inverse to $\theta _{n+1}$, it suffices
  to prove that $\mu _{n+1}\theta _{n+1}(t_\sigma )=t_\sigma $ for $\sigma \in S_{n-1}$.  Let
  $L_\sigma $ denote the closure of $T_\sigma $, which is Milnor's link as
  depicted in Figure \ref{F34}.  Milnor \cite{Milnor} proved that
  for $\tau \in S_{n-1}$
  \begin{equation}
    \mu _\tau (L_\sigma ) = \left\{ \begin{array}{ll}
      1 & \text{if $\tau =\sigma $,} \\
      0 & \text{otherwise.}
    \end{array}\right.
  \end{equation}
  Hence we have
  \begin{equation*}
    \mu _{n+1}\theta _{n+1}(t_\sigma )
    =\sum_{\tau \in S_{n-1}}\bar\mu _\tau (L_\sigma )t_\tau
    =t_\sigma .
  \end{equation*}
  This completes the proof.
\end{proof}

\begin{corollary}
  \label{r37}
  For two Brunnian $(n+1)$-component string links $T,T'\in \BSLnn$, the
  following conditions are equivalent.
  \begin{enumerate}
  \item $T$ and $T'$ are $C^a_{n+1}$-equivalent.
  \item $T$ and $T'$ have the same Milnor invariants of length $n+1$.
  \item $T$ and $T'$ are link-homotopic.
  \end{enumerate}
\end{corollary}

\begin{proof}
  The equivalence $(2)\Leftrightarrow(3)$ is due to Milnor
  \cite{Milnor}.
  The equivalence $(1)\Leftrightarrow(2)$ follows from the proof of
  Proposition \ref{r25}.
\end{proof}

\begin{remark}
  \label{r38}
  Miyazawa and Yasuhara \cite{MY} prove a similar result for Brunnian
  links.  It seems that their proof can be applied to the case of
  string links.  See also the Remark \ref{r42} below.
\end{remark}

\section{The group $\bBnn$}
\label{sec:group-bbnn}

\subsection{The set $\Bnn$}
\label{sec:closing}

Let $\Bnn$ denote the set of the ambient isotopy classes of
$(n+1)$-component Brunnian links.  Let
\begin{equation}
  \label{e14}
  c_{n+1}\zzzcolon \BSLnn\rightarrow \Bnn
\end{equation}
denote the map such that $c_{n+1}(\beta )$ is obtained from $\beta \in \BSLnn$ by
closing each component in the well-known manner.

\begin{proposition}
  \label{r20}
  The map $c_{n+1}$ is onto.
\end{proposition}

\begin{proof}
  This is an immediate consequence of \cite[Proposition 12]{Hb}.
\end{proof}

\subsection{The isomorphism $\bar{c}_{n+1}\zzzcolon \bBSLnn\rightarrow \bBnn$}
\label{sec:isom-bclnnbbslnnbbnn}

Set
\begin{equation*}
  \bBnn=\Bnn/(\text{$C^a_{n+1}$-equivalence}),
\end{equation*}
and let
\begin{equation*}
  \bar{c}_{n+1}\zzzcolon \bBSLnn\rightarrow \bBnn
\end{equation*}
denote the map induced by $c_{n+1}$, which is onto by Proposition
\ref{r20}.

\begin{proposition}
  \label{r39}
  $\bar{c}_{n+1}$ is one-to-one.
\end{proposition}

\begin{proof}
  It suffices to prove that there is a map $\bBnn\rightarrow \modT _{n+1}$ which is
  inverse to $\bar{c}_{n+1}\theta _n\zzzcolon \modT _{n+1}\rightarrow \bBnn$.  This is proved
  similarly as in the proof of Proposition \ref{r25}.
\end{proof}

Proposition \ref{r39} provides the set $\bBnn$ the well-known abelian
group structure, with multiplication induced by band sums of Brunnian
links.

As a corollary, we obtain another proof of a result of Miyazawa and
Yasuhara \cite{MY}.

\begin{corollary}[{\cite[Theorem 1.2]{MY}}]
  \label{r40}
  Let $L$ and $L'$ be two $(n+1)$-component Brunnian links in $S^3$.
  Then the following conditions are equivalent.
  \begin{enumerate}
  \item $L$ and $L'$ are $C^a_{n+1}$-equivalent.
  \item $L$ and $L'$ are $C_{n+1}$-equivalent.
  \item $L$ and $L'$ are link-homotopic.
  \end{enumerate}
\end{corollary}

\begin{proof}
  The result follows immediately from Propositions \ref{r37} and
  \ref{r39}.
\end{proof}

\begin{remark}
  \label{r42}
  Miyazawa and Yasuhara \cite{MY} do not explicitly state the
  equivalence of (1) and others, but this equivalence follows from
  their proof.

  Note that, unlike the $C^a_{n+1}$-equivalence, neither the
  $C_{n+1}$-equivalence nor the link-homotopy are closed for Brunnian
  links.
\end{remark}

\begin{remark}
  \label{r43}
  It is possible to show directly that $\modT _{n+1}$ is isomorphic to
  $\bBnn$, without using string links and the closure map
  $\bar{c}_{n+1}$.  The proof uses Milnor's $\overline{\mu }$-invariants
  and the above result of Miyazawa and Yasuhara.  Our approach
  provides an alternative proof of the latter (instead of using it).
\end{remark}

\subsection{Trees and the Milnor invariants}
\label{sec:trees-miln-invar}

In this subsection, we fix some notations which are used in
later sections.  (Some has appeared in the proof of Proposition
\ref{r25}.)

For $\sigma \in S_{n-1}$, let $t_\sigma $ denote the unitrivalent tree as depicted
in Figure \ref{F32}.  The $t_\sigma $ for $\sigma \in S_{n-1}$ form a basis of
$\modT _{n+1}$.  Let $T_\sigma $ denote the corresponding $C^a_n$-tree for the
$(n+1)$-component unlink $U=U_1\cup \dots \cup U_{n+1}$, see Figure \ref{F33}.
{\def\te#1{$U_{\sigma (#1)}$}\FI{F33}{}}

For $i_1,\ldots ,i_{n+1}$ with $\{i_1,\ldots ,i_{n+1}\}=\{1,\ldots ,n+1\}$, let
\begin{equation*}
  \bar\mu _{i_1,\ldots ,i_{n+1}}\zzzcolon \Bnn\rightarrow \modZ
\end{equation*}
denote the Milnor invariant, which is additive under connected sum
\cite{Milnor} (see also \cite{Cochran,Orr,Krushkal}).  For $\sigma \in S_{n-1}$,
we set
\begin{equation*}
  \bar\mu _\sigma = \bar\mu _{\sigma (1),\sigma (2),\ldots ,\sigma (n-1),n,n+1}\zzzcolon \Bnn\rightarrow \modZ .
\end{equation*}
It is well known \cite{Milnor} that for $\rho \in S_{n-1}$
\begin{equation*}
  \bar\mu _\rho (U_{T_\sigma }) = \left\{ \begin{array}{ll}
    1 & \text{if $\rho =\sigma $,}\\
    0 & \text{otherwise.}
  \end{array}\right.
\end{equation*}

\section{The Goussarov-Vassiliev filtration for links}
\label{sec:gouss-vass-filtr}

In this section, we briefly recall the formulation using claspers of
the Goussarov-Vassiliev filtrations for links.  See \cite{H}
for details.

\subsection{Forest schemes and Goussarov-Vassiliev filtration}
\label{sec:forest-schem-gouss}

A {\em forest scheme} of degree $k$ for a link $L$ in a $3$-manifold
$M$ will mean a collection $S=\{G_1,\ldots ,G_l\}$ of disjoint
(strict) tree claspers $G_1,\ldots ,G_l$ for $L$ such that
$\sum_{i=1}^k\deg G_i=k$.  A forest scheme $S$ is said to be {\em
simple} if every element of $S$ is simple.

For $n\ge 0$, let $\modL(M,n)$ denote the set of ambient isotopy classes
of oriented, ordered links in $M$.

For a forest scheme $S=\{G_1,\ldots ,G_l\}$ for a link $L$ in $M$, we
set
\begin{equation*}
  [L,S] = [L;G_1,\ldots ,G_l]
  = \sum_{S'\subset S} (-1)^{|S'|} L_{\bigcup S'} \in \modZ \modL (M,n),
\end{equation*}
where the sum is over all subsets $S'$ of $S$, and $|S'|$ denote the
number of elements of $S'$.

For $k\ge 0$, let $J_k(M,n)$ (sometimes denoted simply by $J_k$) denote
the $\modZ $-submodule of $\modZ\modL(M,n)$ generated by the elements of
the form $[L,S]$, where $L\in \modL (M,n)$ and $S$ is a forest scheme
for $L$ of degree $k$.  We have
\begin{equation*}
  \modZ\modL(M,n)=J_0(M,n)\subset J_1(M,n)\subset \cdots,
\end{equation*}
which coincides with the Goussarov-Vassiliev filtration using
alternating sums of links determined by singular links, see
\cite[Section 6]{H}.

\subsection{Crossed edge notation}
\label{sec:cross-edge-notat}

It is useful to introduce a notation for depicting certain linear
combinations of surgery along claspers, which we call {\em crossed
edge notation}.

Let $G$ be a clasper for a link $L$ in a $3$-manifold $M$.  Let $E$ be
an edge of $G$.  By putting a cross on the edge $E$ in a figure, we
mean the difference $L_G-L_{G_0}$, where $G_0$ is obtained from $G$ by
inserting two trivial leaves into $E$.  See Figure \ref{F29}.
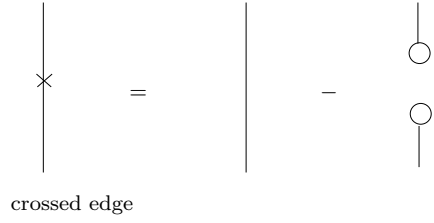
\FI{F29}{The crossed edge notation} If we put several crosses on the
edges of $G$, then we understand it in a multilinear way.  I.e., a
clasper with several crosses is an alternating sum of the result of
surgery along claspers obtained from $G$ by inserting pairs of
trivial, unlinked leaves into the crossed edges.  We will freely use
the identities depicted in Figure \ref{F30}, which can be easily
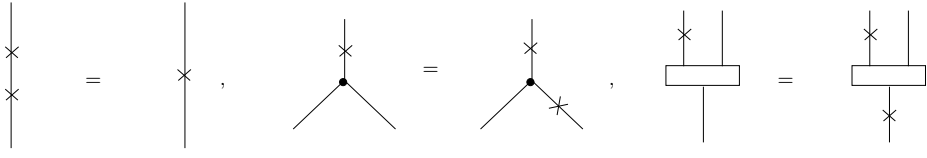
verified.  \FI{F30}{Identities for the crossed edge notation} The
second identity implies that if $G'$ is a connected graph clasper
contained in $G$ and there are several crosses on $G'$, then one can
safely replace these crosses by just one cross on one edge in $G'$.
This properties can be generalized to the case where $G'$ is a
connected subsurface of $G$ consisting only of nodes, edges, leaves
and disk-leaves.  Note also that if $S=\{G_1,\ldots ,G_l\}$, is a forest
scheme for $L$, then $[L,S]$ can be expressed by the clasper
$G_1\cup \dots \cup G_l$ with one cross on each component $G_i$.

\section{Goussarov-Vassiliev invariants of Brunnian links}
\label{sec:gouss-vass-invar}

Throughout this section, let $U=U_1\cup U_2\cup \cdots\cup U_{n+1}$ be
the $(n+1)$-component unlink.

\subsection{The map $\kappa _{n+1}\zzzcolon \bBnn\rightarrow \bJ_{2n}(n+1)$}
\label{sec:map-n+1}

\begin{proposition}
  \label{r27}
  Let $n\ge 2$.  Let $L$ and $L'$ be two $(n+1)$-component Brunnian
  links in an oriented, connected $3$-manifold $M$.  If $L$ and $L'$
  are $C^a_{n+1}$-equivalent (or link-homotopic), then we have
  $L'-L\in J_{2n+1}$.
\end{proposition}

Proposition \ref{r27} implies the following.

\begin{corollary}
  \label{r30}
  The restriction of any Goussarov-Vassiliev invariant
  of degree $2n$ to $(n+1)$-component Brunnian links is a
  link-homotopy invariant.
\end{corollary}

\begin{proof}[Proof of Proposition \ref{r27}]
  First, we consider the case $L=U$.  By using the same arguments as
  in the proof of \cite[Lemma 14]{Hb}, we see that there is a clasper
  $G$ for $U$ consisting of $C^a_l$-claspers with $n+1\le l<2n+1$,
  such that $U$ bounds $n+1$ disjoint disks which are disjoint from
  the edges and the nodes of $G$, and such that $U_G\sim_{C^a_{2n+1}}
  U_T$.  The latter implies that $U_G - U_T\in J_{2n+1}$.  We use the
  equality $U_G=\sum_{G'\subset G} [U,G']$.  Clearly
  $[U,G']\in J_{2n+1}$ for $|G'|>1$, so we may safely assume that $G$
  has only one component.  We then have $U_G - U\in J_{2n+1}$ as a
  direct application of \cite[Lemma 16]{Hb}.  This completes the proof
  of the case $L=U$.

  Now consider the general case.  We may assume that $L'$ is obtained
  from $L$ by one simple $C^a_n$-move.  Since $L$ is an
  $(n+1)$-component Brunnian link, it follows from Theorem \ref{r12}
  and Lemma \ref{r18} that there exists a simple $C^a_n$-forest $F$
  for $U$ such that $L=U_F$.  Also, there exists a simple
  $C^a_{n+1}$-tree $\tilde{T}$ for $L=U$ such that $L'=L_{\tilde T}$.
  We may assume that $\tilde T$ is a simple $C^a_{n+1}$-tree for $U$
  disjoint from $F$ such that $L'=U_{F\cup \tilde{T}}$.  Let $S$ be the
  forest scheme consisting of the trees $T_1,\ldots,T_l$ of $F$.  We
  have $L=\sum_{S'\subset S} [U,S']$ and $L'= \sum_{S'\subset S}
  [U_T,S']$.  Hence we have
  \begin{equation*}
    \begin{split}
      L' - L
      = \sum_{S'\subset S} [U,S'\cup \{T\}].
    \end{split}
  \end{equation*}
  Since $\deg \tilde T=n+1$ and $\deg T_i=n$ for all $i$, the term in
  the above sum is contained in $J_{2n+1}$ unless $S'=\emptyset$.
  Hence we have
  \begin{equation*}
    L'-L\equiv [U,\tilde T]\equiv0 \pmod{J_{2n+1}},
  \end{equation*}
  where the second congruence follows from the first case.
\end{proof}

By Proposition \ref{r27}, we have a map
\begin{equation*}
   \kappa _{n+1}: \bBnn \longrightarrow \bJ_{2n}(n+1)
\end{equation*}
defined by $\kappa _{n+1}(L)=[L-U]_{J_{2n+1}}$.

\subsection{Quadraticity of $\kappa _{n+1}$}
\label{sec:quadraticity-n+1}

Let $n\ge 2$.  In this subsection, we establish the following commutative
diagram.
\begin{equation}
  \label{e5}
  \begin{CD}
    \modT _{n+1} @>\psi _{n+1}>\simeq> \bBnn \\
    @Vq_{n+1}VV @VV\kappa _{n+1}V\\
    \SymT @>>\delta _{n+1}> \bJ_{2n}(n+1)
  \end{CD}
\end{equation}
Definitions of $\psi _{n+1}$, $\SymT$, $q_{n+1}$ and $\delta _{n+1}$ are
in order.

The isomorphism $\psi _{n+1}$ is the composition of
\begin{equation*}
  \modT _{n+1}\xyto{\theta _n}{\simeq}\bBSLnn\xyto{\bar{c}_{n+1}}{\simeq}\bBnn.
\end{equation*}

Let $\Sym^2\modT ^\modQ _{n+1}$ denote the symmetric product of two
copies of $\modT ^\modQ _{n+1}:=\modT _{n+1}\otimes \modQ $, and let
$\SymT$ denote the $\modZ $-submodule of $\Sym^2\modT ^\modQ _{n+1}$
generated by $\half x^2$, $x\in \modT _{n+1}$.  One can easily verify
that $\SymT$ is $\modZ $-spanned by the elements $\half t_\sigma ^2$
for $\sigma \in S_{n-1}$ and $t_\sigma t_{\sigma '}$ for $\sigma
,\sigma '\in S_{n-1}$.  (Of course we have $t_\sigma t_{\sigma
'}=t_{\sigma '}t_\sigma $.  Thus $\SymT$ is a free abelian group of
rank $\half(n-1)!((n-1)!+1)$.)

The arrow $q_{n+1}$ is the quadratic map defined by $q_{n+1}(x)=\half
x^2$ for $x\in \modT _{n+1}$.

The arrow $\delta _{n+1}$ is the homomorphism defined as follows.  For
$\sigma ,\sigma '\in S_{n-1}$, let $T_\sigma $ and $T_{\sigma '}$ be
the corresponding simple $C^a_n$-trees for $U$ as in Section
\ref{sec:trees-miln-invar}.  Let $\tilde T_{\sigma '}$ denote a simple
$C^a_n$-trees obtained from $T_{\sigma '}$ by a small isotopy if
necessary so that $\tilde T_{\sigma '}$ is disjoint from $T_\sigma $.
Set
\begin{equation*}
  \delta _{n+1}(t_\sigma t_{\sigma '})
  =[U;T_\sigma ,\tilde T_{\sigma '}]_{J_{2n+1}}\in \bJ_{2n}(n+1),
\end{equation*}
which does not depend on how we obtained $\tilde T_{\sigma '}$ from
$T_\sigma $, since crossing changes between an edge of $T_\sigma $ and
an edge of $\tilde T_{\sigma '}$ preserves the right-hand side.  (This
can be verified by using a `$C^a$-version' of \cite[Proposition
4.6]{H}.)  For the case of $\half t_\sigma ^2$, we modify the above
definition with $\sigma '=\sigma $ as follows.  Let $T_\sigma $ and
$\tilde T_{\sigma }$ be as above.  See Figure \ref{F17}.  \FI{F17}{}
Let $T'$ be the $C_{n-1}$-tree obtained from $\tilde T_\sigma $ by
first removing the disk-leaf $D$ intersecting $U_{n+1}$, the edge $E$
incident to $D$, and the node $N$ incident to $E$, and then gluing the
ends of the two edges which were attached to $N$.  Moreover, let $C$
be a $C_1$-tree which intersects $U_{n+1}$ and $U_{\sigma (1)}$ as
depicted.  Set
\begin{equation*}
  \delta _{n+1}(\half t_\sigma ^2)= [U;T_\sigma ,T',C].
\end{equation*}

\begin{lemma}
  \label{r3}
  We have
  \begin{equation}
    \label{e6}
    [U;T_\sigma ,\tilde T_\sigma ]\equiv 2[U;T_\sigma ,T',C]   \pmod{J_{2n+1}}.
  \end{equation}
\end{lemma}

\begin{proof}
  By \cite[Section 8.2]{H}, it suffices to prove the identity in the
  space of unitrivalent diagram depicted in Figure \ref{F28}, which
  can be easily verified using the STU relation several times.
  \FI{F28}{}
\end{proof}

It follows from Lemma \ref{r3} that $\delta _{n+1}$ is a well-defined
homomorphism.  Set
\begin{equation*}
  \half [U;T_\sigma ,\tilde T_\sigma ]_{J_{2n+1}}=[U;T_\sigma ,T',C]_{J_{2n+1}}.
\end{equation*}
We have
\begin{equation*}
  \delta _{n+1}(\half t_\sigma ^2)= \half [U;T_\sigma ,\tilde T_\sigma ]_{J_{2n+1}}.
\end{equation*}

\begin{theorem}
  \label{r15}
  The diagram \eqref{e5} commutes.  In particular, $\kappa _{n+1}$ is a
  quadratic map.
\end{theorem}

We need the following lemma before proving Theorem \ref{r15}.

\begin{lemma}
  \label{r5}
  Let $C$ be a clasper for a link $L$ such that there is a disk-leaf
  $D$ of $T$ which `monopolizes' a component $K$ of $L$ in the sense
  of \cite[Definition 15]{Hb}, and such that $D$ is adjacent to a
  node.  That is, $T$ and $L$ looks as depicted in the left hand side
  of Figure \ref{F19}.  \FI{F19}{}  Then we have the identity as
  depicted in the figure.
\end{lemma}

\begin{proof}
  The identity is easily verified and left to the reader.  (Note that
  Lemma \ref{r5} is essentially the same as \cite[(4.4)]{Hb}.)
\end{proof}

\begin{proof}[Proof of Theorem \ref{r15}]
  Let $\sigma \in S_{n-1}$.  We must show that
  \begin{equation*}
    [U;T_\sigma ]_{J_{2n+1}}= \half [U;T_\sigma ,\tilde T_\sigma ]_{J_{2n+1}}.
  \end{equation*}

  For $i=1,\ldots ,n+1$, let $D_i$ denote the disk-leaf of $T_\sigma $
  intersecting $L_i$, and let $E_i$ denote the incident edge.
  For $i=1,\ldots ,n-1$, let $N_i$ denote the node incident to $E_i$.

  By applying Lemma \ref{r5} to the edge of $T_\sigma $ which is incident
  to $N_{\sigma (1)}$ but not to $D_{n+1}$ or $D_{\sigma (1)}$, we obtain the
  identity depicted in Figure \ref{F18}.  \FI{F18}{} Let $B$ be the
  box and $E$ be the edge as depicted.  Let $G$ be the clasper in the
  right hand side.  By zip construction \cite[Section 3.3]{H} at $E$,
  we obtain a crossed clasper depicted in Figure \ref{F20new}, which
  consists of two components $T_\sigma $ and $P$.  \FI{F20new}{} The component
  $P$ has $n-2$ (non-disk) leaves.

  We claim that we can unlink the leaves of $P$ from $T_\sigma $
  without changing the class in $\bJ_{2n}(n+1)$.  To see this, it
  suffices to show that
  \begin{equation}
    \label{e7}
    U_{T_\sigma \cup P}-U_{T_\sigma \cup P'}\in J_{2n+1},
  \end{equation}
  where $P'$ is obtained from $P$ by the unlinking operation.  Note
  that each unlinking is performed by a sequence of crossing changes
  between an edge of the $C^a_n$-tree $T_\sigma $ and a link component
  (after performing surgery along $P$ in the regular neighborhood of
  $P$), and thus can be performed by $C^a_{n+1}$-moves.  Since all the
  links appearing in this sequence is Brunnian, we have \eqref{e7} by
  Proposition \ref{r27}.  This completes the proof of the claim.

  By the above claim, it follows that
  \begin{equation*}
  [U;T_\sigma ,P]\equiv[U;T_\sigma ,T']\pmod{J_{2n+1}},
  \end{equation*}
  where $T'$ is obtained from $P$ by removing the leaves, the incident
  edges, and the boxes, and then smoothing the open edges, see the
  left hand side of Figure \ref{F21new}, which is equal to the right hand
  side by Lemma \ref{r5}. \FI{F21new}{} The result is related to the
  desired clasper defining $[U;T_\sigma ,T',C]$ by half twists of two
  edges and homotopy with respect to $U$, and hence equivalent modulo
  $J_{2n+1}$ to $[U;T_\sigma ,T',C]$.  This completes the proof.
\end{proof}

\subsection{Proof of Theorem \ref{r6}}
\label{sec:proof-theorem-refr6}

In this subsection we prove Theorem \ref{r6}.

Let $L\in \Bnn$.  We have
\begin{equation*}
  [L]_{C^a_{n+1}} = \sums \bar\mu _\sigma (L) [U_{T_\sigma }]_{C^a_{n+1}}
\end{equation*}
in $\bBnn$.  (Recall that the sum is induced by band-sum in
$\bBnn$.)  Hence we have by the commutativity of \eqref{e5}
\begin{equation*}
  \begin{split}
    [L-U]_{J_{2n+1}}
    &= \kappa _{n+1}([L]_{C^a_{n+1}})\\
    &= \delta _{n+1}q_{n+1}\psi _{n+1}^{-1}
    (\sums \bar\mu _\sigma (L) [U_{T_\sigma }]_{C^a_{n+1}})\\
    &= \delta _{n+1}q_{n+1}(\sums \bar\mu _\sigma (L) t_\sigma )\\
    &= \delta _{n+1}(\half(\sums \bar\mu _\sigma (L) t_\sigma )^2)\\
    &= \delta _{n+1}(\half\sumss \bar\mu _\sigma (L)\bar\mu _{\sigma '}(L) t_\sigma t_{\sigma '})\\
    &= \half\sumss \bar\mu _\sigma (L)\bar\mu _{\sigma '}(L) [U;T_\sigma ,\tilde T_{\sigma '}]_{J_{2n+1}}.
  \end{split}
\end{equation*}
Hence we have
\begin{equation}
  \label{e9}
  f(L)-f(U)=\half\sumss \bar\mu _\sigma (L)\bar\mu _{\sigma '}(L) f([U;T_\sigma ,\tilde T_{\sigma '}]).
\end{equation}
We give any total order on the set $S_{n-1}$.  Then we have
\begin{equation*}
  \begin{split}
  f(L)-f(U)&=\sums(\half f([U;T_\sigma ,\tilde T_\sigma ]))\bar\mu _\sigma (L)\bar\mu _\sigma (L)\\
  &\quad \quad +
  \sum_{\sigma <\sigma '} f([U;T_\sigma ,\tilde T_{\sigma '}]) \bar\mu _\sigma (L)\bar\mu _{\sigma '}(L) .
  \end{split}
\end{equation*}
Note that $\half f([U;T_\sigma ,\tilde T_\sigma ])\in \modZ $ and $f([U;T_\sigma ,\tilde T_{\sigma '}])\in \modZ $.
Hence we have \eqref{e8} by setting
\begin{equation*}
  f_{\sigma ,\sigma '} =
  \begin{cases}
    \half f([U;T_\sigma ,\tilde T_\sigma ]) & \text{if $\sigma =\sigma '$},\\
    f([U;T_\sigma ,\tilde T_{\sigma '}]) & \text{if $\sigma <\sigma '$},\\
    0 & \text{if $\sigma >\sigma '$}.
  \end{cases}
\end{equation*}

This completes the proof of Theorem \ref{r6}.

\end{document}